\documentclass{amsart}
\usepackage{amsmath, amsthm, amssymb, amscd, epsfig}

\begin{document}

\newtheorem{theorem}{Theorem}[section]
\newtheorem{lemma}[theorem]{Lemma}
\newtheorem{proposition}[theorem]{Proposition}
\newtheorem{corollary}[theorem]{Corollary}
\newtheorem{conjecture}[theorem]{Conjecture}
\newtheorem{question}[theorem]{Question}
\newtheorem{problem}[theorem]{Problem}
\newtheorem*{claim}{Claim}
\newtheorem*{criterion}{Criterion}
\newtheorem*{basmajian}{Basmajian's Identity}
\newtheorem*{bridgeman}{Bridgeman's Identity}
\newtheorem*{chimney}{Chimney Decomposition}
\newtheorem*{rational_thm}{Rational Identity}

\theoremstyle{definition}
\newtheorem*{definition}{Definition}
\newtheorem{construction}[theorem]{Construction}
\newtheorem{notation}[theorem]{Notation}

\theoremstyle{remark}
\newtheorem*{remark}{Remark}
\newtheorem*{example}{Example}

\def\Z{\mathbb Z}
\def\H{\mathbb H}

\def\area{\textnormal{area}}
\def\volume{\textnormal{volume}}

\def\til{\tilde}

\title{Chimneys, leopard spots, and the identities of Basmajian and Bridgeman}
\author{Danny Calegari}
\address{Department of Mathematics \\ Caltech \\
Pasadena CA, 91125}
\date{\today}

\begin{abstract}
We give a simple geometric argument to derive in a common manner orthospectrum identities of
Basmajian and Bridgeman. Our method also considerably simplifies the determination of the
summands in these identities. For example, for every odd integer $n$, there is a rational function
$q_n$ of degree $2(n-2)$ so that if $M$ is a compact hyperbolic manifold of dimension $n$ with
totally geodesic boundary $S$, there is an identity $\chi(S) = \sum_i q_n(e^{l_i})$ where the
sum is taken over the orthospectrum of $M$. When $n=3$, this has the explicit form 
$\sum_i 1/(e^{2l_i}-1) = -\chi(S)/4$.
\end{abstract}

\maketitle

\section{Orthospectrum identities}

Let $M$ be a compact hyperbolic $n$-manifold with totally geodesic boundary $S$. An {\em orthogeodesic}
is a properly immersed geodesic arc perpendicular to $S$ at either end. The {\em orthospectrum} is the
set of lengths of orthogeodesics, counted with multiplicity.

\smallskip

Basmajian \cite{Basmajian} and Bridgeman--Kahn \cite{Bridgeman, Bridgeman_Kahn} derived identities relating the orthospectrum
of $M$ to the area of $S$ and the volume of $M$ respectively. The following identity
is implicit in \cite{Basmajian}:

\begin{basmajian}[\cite{Basmajian}]
There is a function $a_n$ depending only on $n$, so that if $M$ is a compact hyperbolic $n$-manifold
with totally geodesic boundary $S$, and $l_i$ denotes the (ordered) orthospectrum of $M$, with multiplicity,
there is an identity
$$\area(S) = \sum_i a_n(l_i)$$
\end{basmajian}

Basmajian's identity is not well known; in fact, Bridgeman and Kahn were apparently unaware of 
Basmajian's work when they derived the following by an entirely different method:

\begin{bridgeman}[\cite{Bridgeman, Bridgeman_Kahn}]
There is a function $v_n$ depending only on $n$, so that if $M$ is a compact hyperbolic $n$-manifold
with totally geodesic boundary $S$, and $l_i$ denotes the (ordered) orthospectrum of $M$, with multiplicity,
there is an identity
$$\volume(M) = \sum_i v_n(l_i)$$
\end{bridgeman}

In this paper, we show that both theorems can be derived from a common geometric perspective. In
fact, the derivation gives a very simple expression for the functions $a_n$ and $v_n$, which we describe
in \S~\ref{rationality_section}. The derivation rests on a simple geometric decomposition.

\begin{definition}
Let $\pi$ and $\pi'$ be totally geodesic $\H^{n-1}$'s in $\H^n$ with disjoint closure in $\H^n \cup S^{n-1}_\infty$. 
A {\em chimney} is the closure of the union of the geodesic arcs from $\pi$ to $\pi'$ 
that are perpendicular to $\pi$. 
\end{definition}

Thus, the boundary of the chimney consists of three pieces: the {\em base}, which is a round disk in $\pi$,
the {\em side}, which is a cylinder foliated by geodesic rays, and the {\em top}, which is the plane $\pi'$.
Note that the distance from the base to the top is realized by a unique orthogeodesic, called the {\em core}. 
The {\em height} of the chimney is the length of this orthogeodesic, and the {\em radius} is the 
radius of the base (these two quantities are related, and either one determines the chimney up to isometry).

\begin{chimney}
Let $M$ be a compact hyperbolic $n$-manifold with totally geodesic boundary $S$. Let $M_S$ be the covering
space of $M$ associated to $S$. Then $M_S$ has a canonical decomposition into a piece of zero measure,
together with two chimneys of height $l_i$ for each number $l_i$ in the orthospectrum.
\end{chimney}
\begin{proof}
If $S$ is disconnected, the cover $M_S$ is also disconnected, and consists of a union of connected covering
spaces of $M$, one for each component of $S$. The boundary of $M_S$ consists of a copy of $S$, together
with a union of totally geodesic planes. Each such plane is the top of a chimney, with base a round disk in $S$,
and these chimneys are pairwise disjoint and embedded. Since $M$ is geometrically finite, the limit 
set has measure zero, and therefore these chimneys exhaust all of $M_S$ except for a subset of measure zero.
Every oriented orthogeodesic in $M$ lifts to a unique geodesic arc with initial point in $M_S$. 
Evidently this arc is the core of a unique chimney in the decomposition, and all chimneys arise this way.
\end{proof}

Basmajian's identity is immediate (in fact, though Basmajian does not express things in these terms,
the argument we give is quite similar to his):

\begin{proof}
$S$ in $M_S$ is decomposed into a set
of measure zero together with the union of the bases of the chimneys. Thus
$$\area(S) = 2\sum_i \text{ area of the base of a chimney of height }l_i$$
\end{proof}

\begin{remark}
Thurston calls the chimney bases {\em leopard spots}; they arise in the definition of the skinning map
(see e.g. \cite{Otal}).
\end{remark}

Bridgeman's identity takes slightly more work, but is still elementary:

\begin{proof}
If $p$ is a point in $M$, and $\gamma$ is an arc from $p$ to $S$, there is a unique geodesic in the
relative homotopy class of $p$ which is perpendicular to $S$. Thus, the unit tangent sphere to $p$ is
decomposed into a set of measure zero, together with a union of round disks, one for each relative
homotopy class of arc $\gamma$. 

The area of the disk in $UT_p$ associated to $\gamma$ can be computed as
follows. Let $\til{\gamma}$ be the unique lift of $\gamma$ to $M_S$ with one endpoint on $S$, and let
$\til{p}$, a lift of $p$, be the other endpoint of $\til{\gamma}$. If $N$ is the complete hyperbolic manifold with
$M$ as compact core and $N_S$ denotes the cover of $N$ 
associated to $S$ (so that $M_S$ is a convex subset of $N_S$), 
let $h_S$ be the harmonic function on $N_S$ whose value at every point $q$
is the probability that Brownian motion starting at $q$ exits the end associated to $S$. Note that
$h_S=1/2$ on $S$, and at every point $q$ depends only on the distance from $q$ to $S$. Then the
area of the disk in $UT_p$ associated to $\gamma$ is $\Omega_{n-1}\cdot h_S(\til{p})$, where 
$\Omega_{n-1}:=2\pi^{n/2}/\Gamma(n/2)$ denotes the area of a Euclidean sphere of dimension $n-1$ 
and radius $1$.

Since the volume of the unit tangent bundle of $M$ is $\Omega_{n-1}\cdot\volume(M)$, 
it follows that the volume of $M$ is equal to the integral of $h_S$ over $M_S$. In each chimney,
$h_S$ restricts to a harmonic function $h$, equal to $1/2$ on the base, and whose value at each point
depends only on the distance to the base. Hence
$$\volume(M) = 2\sum_i \text{ integral of }h\text{ over a chimney of height }l_i$$
\end{proof}

\begin{remark}
In fact, precisely because our derivation is utterly unlike that of \cite{Bridgeman_Kahn}, we do
{\em not} know whether Bridgeman's function $v_n$ is equal to the integral of $h$ over an $n$-dimensional
chimney of given height, only that there {\em is} such a function $v_n$ with the desired properties. 
If $n=2$, our $v_2$ and Bridgeman's $v_2$ agree, but the proof is not easy; one short
derivation follows from \cite{Calegari_ortho}, together with a geometric dissection argument.
\end{remark}

\section{Explicit formulae}\label{rationality_section}

In this section we show that the summands in the area and volume identities have a very nice explicit
form. The expressions we obtain depend on the following elementary ingredients:

\smallskip

{\noindent \bf quadrilateral}: A chimney is a solid of revolution, obtained by revolving a hyperbolic 
quadrilateral $Q$ with three right angles and one ideal vertex about the $S^{n-2}$ of directions
perpendicular to one of the finite sides (which becomes the core of the chimney, the other finite
side becoming the radius of the base). In a quadrilateral with three right angles and one ideal
vertex, the length of one finite edge determines the other. If one finite edge has length $l$, 
let $\iota(l)$ denote the length of the other finite edge, so that $\iota$ is an involution on
$(0,\infty)$. Then $\iota$ is defined implicitly by 
the fact that it is positive, and the identity
$$1/\cosh^2(l) + 1/\cosh^2(\iota(l)) = 1$$
or equivalently,
$$\sinh(\iota(l)) = 1/\sinh(l)$$
If we write $\alpha = e^l$ and $\beta = e^{\iota(l)}$, then $\alpha$ and $\beta$ are related by
$$\beta + \beta^{-1} = 2\left( \frac {\alpha + \alpha^{-1}} {\alpha - \alpha^{-1}} \right)$$

\smallskip

{\noindent \bf hyperbolic volume:} If $B$ is a ball of radius $r$ in $n$-dimensional hyperbolic space, 
let $V_n^H(r)$ denote the volume of $B$. One has the following integral formula for $V_n^H$:
$$V_n^H(r) = \Omega_{n-1} \int_0^r \sinh^{n-1}(t) dt$$

\smallskip

The base of an $n$-dimensional chimney of height $l$ is just the volume of an $(n-1)$-dimensional
ball in hyperbolic space of radius $\iota(l)$.
For $n$ even, the integral $\int_0^{\iota(l)} \sinh^{n-1}(t) dt$ is a 
{\em polynomial} in $\beta+\beta^{-1}$, and therefore a {\em rational
function} in $\alpha$ of degree $2(n-1)$. If the dimension of $M$ is at least $3$, the set of numbers
$e^l$ where $l$ runs over the orthospectrum are algebraic (by Mostow rigidity), 
and contained in a quadratic extension of the trace field of $M$.

If $S$ has even dimension, then the area of $S$ is proportional to the Euler characteristic, by the
Chern--Gauss--Bonnet theorem; in fact, for a hyperbolic manifold of dimension $n$ where $n$ is even,
one has:
$$\area(S) = (2\pi)^{n/2}\chi(S)r_n$$
where $r_n$ is a rational number depending on $n$.

The following corollary appears to be new:

\begin{rational_thm}
For every odd integer $n$, there is a rational function $q_n$ of degree $2(n-2)$, 
with integral coefficients, so that if $M$ is a compact hyperbolic manifold of (odd) 
dimension $n$ with totally geodesic boundary $S$, there is an identity
$$\chi(S) = \sum_i q_n(e^{l_i})$$
where $\chi$ denotes Euler characteristic (which takes values in $\Z$) and
$l_i$ denotes the orthospectrum of $M$ (with multiplicity). Note that for $n \ge 3$,
the numbers $e^{l_i}$ are all contained in a fixed number field $K$ (depending on $M$).
\end{rational_thm}

\begin{example}
It is elementary to compute $q_n$ for small $n$. For example,
$$q_3(x) = \frac {4} {1-x^2}$$
$$q_5(x) = \frac {5x^6 - 33x^4 + 63x^2-27} {8(x^2 - 1)^3}$$
The denominator is easily seen to be an integer multiple of $(x^2-1)^{n-2}$.
\end{example}

\begin{remark}
In the case of $3$ dimensions, the identity has the following form. Let $M$ be 
a hyperbolic $3$-manifold with totally geodesic boundary $S$. Then
$$\sum_i \frac 1 {e^{2l_i} -1 } = -\chi(S)/4$$ 
This is vaguely reminiscent of McShane's identity \cite{McShane}, which says that for $S$ a hyperbolic
once-punctured torus, there is an identity
$$\sum_i \frac 1 {1+e^{l_i}} = 1/2$$
where the sum is taken over lengths $l_i$ of {\em simple} closed geodesics {\em in the surface $S$}.

If there is a simple relation between our identities and McShane's identity, it is not obvious.
However, Mirzakhani \cite{Mirzak} showed how to derive and generalize McShane's identity as a
sum over {\em embedded} orthogeodesics on a surface with boundary. The appearance of orthogeodesics
in yet another identity is quite suggestive of a more substantial connection, though we do not 
know what it might be.
\end{remark}

To determine the summands in the volume identity, one needs the following additional ingredients:

\smallskip

{\noindent \bf $\phi$-quadrilateral:} If $Q$ is a hyperbolic quadrilateral with three right angles and one
vertex with angle $\phi$, then one of the lengths $l$ of the edges ending at right angles
determines the other $\iota_\phi(l)$, defined implicitly by the identity
$$\sinh(\iota_\phi(l)) = \sinh(\iota(l))\cos(\phi) = \cos(\phi)/\sinh(l)$$

\smallskip

{\noindent \bf spherical volume:} If $B$ is a ball of radius $r$ in $n$-dimensional spherical space, 
let $V_n^S(r)$ denote the volume of $B$. One has the following integral formula for $V_n^S$:
$$V_n^S(r) = \Omega_{n-1} \int_0^r \sin^{n-1}(t) dt$$

\smallskip

{\noindent \bf harmonic:} Let $h$ be the harmonic function on $\H^n$ equal to the indicator function
of a round disk $D$ in $S^{n-1}_\infty$, so that $h=1/2$ on the plane $\pi$ bounded by $\partial D$.
For $q$ bounded away from $D$ by $\pi$, if $t$ is the distance from $q$ to $\pi$, then $h(q)$ is
$\Omega_{n-1}^{-1}$ times the volume of a ball in $S^{n-1}$ of radius $\theta$, 
where $\sin(\theta) = 1/\cosh(t)$.

\smallskip

{\noindent \bf level sets:} Nearest point projection from an equidistant surface to a totally geodesic
hyperplane multiplies distances by $1/\cosh(t)$. If $C$ is a chimney of height $l$ (and radius $\iota(l)$), 
let $C_t$ be the level set at distance $t$ from the base. Orthogonal projection of $C_t$ to the base of the
chimney is surjective if $t \le l$, and otherwise surjects onto an annulus with outer radius $\iota(l)$,
and inner radius $\iota_\phi(l)$, where $\phi$ is defined implicitly by $\sin(\phi) = \cosh(l)/\cosh(t)$.

The area of $C_t$ is therefore
$$\area(C_t) = \begin{cases}
\cosh^{n-1}(t)V_{n-1}^H(\iota(l)) & \text{ if } t\le l \\
\cosh^{n-1}(t)(V_{n-1}^H(\iota(l)) - V_{n-1}^H(\iota_\phi(l) )) & \text{ if } t\ge l \\
\end{cases}$$

\smallskip

Putting this all together, we get an explicit integral formula for $v_n$ :

\begin{align*}
v_n(l)/2 & = \int_0^l \cosh^{n-1}(t)V_{n-1}^H(\iota(l))V_{n-1}^S(\arcsin(1/\cosh(t)))\Omega_{n-1}^{-1} dt\\
&+\int_l^\infty \cosh^{n-1}(t)(V_{n-1}^H(\iota(l)) - V_{n-1}^H(\iota_\phi(l) ))V_{n-1}^S(\arcsin(1/\cosh(t)))\Omega_{n-1}^{-1} dt\\
\end{align*}

Notice when $n$ is even this can be expressed in closed form in terms of elementary functions
(compare with the formulae and the derivation in \cite{Bridgeman_Kahn}, pp. 4--11).

\section{Acknowledgments}

I would like to thank Martin Bridgeman, Jeremy Kahn, Sadayoshi Kojima, Greg McShane and 
Maryam Mirzakhani for some useful discussions. In particular, this paper owes an obvious debt to
\cite{Bridgeman} and the beautiful sequel \cite{Bridgeman_Kahn}. Thanks also to the referee
for a useful correction. The first version of this paper was written before the author was 
aware of \cite{Basmajian}, and I am very grateful to Greg and Sadayoshi for bringing it to my attention.
Danny Calegari was supported by NSF grant DMS 0707130.

\end{document}